\documentclass[11pt, letterpaper]{article}
\usepackage{authblk}
\usepackage{mathtools} 
\usepackage{amsmath, amssymb} 
\usepackage{amsthm}
\usepackage{enumerate}  
\usepackage{hyperref}
\usepackage{url} 

\usepackage{tikz}
\usepackage{tkz-graph}
\usetikzlibrary{shapes}
\usetikzlibrary{arrows}
\usetikzlibrary{decorations.markings}

\usepackage{graphicx}
\usepackage{caption,subcaption}

\usepackage{geometry}
\newcommand\cx{{\mathbb C}}

\newcommand\ints{{\mathbb Z}}
\newcommand\re{{\mathbb R}}
\newcommand\rats{{\mathbb Q}}

\DeclarePairedDelimiter\abs{\lvert}{\rvert}%
\DeclarePairedDelimiter\norm{\lVert}{\rVert}%

\makeatletter
\let\oldabs\abs
\def\abs{\@ifstar{\oldabs}{\oldabs*}}
\let\oldnorm\norm
\def\norm{\@ifstar{\oldnorm}{\oldnorm*}}
\makeatother

\newcommand\opk[1]{\mathop{\mathrm{#1}}\nolimits}

\newcommand\sbs{\subseteq}

\newcommand\comp[1]{{\mkern2mu\overline{\mkern-2mu#1}}}
\newcommand\pmat[1]{\begin{pmatrix} #1 \end{pmatrix}}
\newcommand\seq[4]{#1_{#2},#1_{#3},\ldots,#1_{#4}}

\newtheoremstyle{plainsl}%
	{\topsep}
	{\topsep}
	{\slshape} 
	{}
	{\normalfont\bfseries}
	{.}
	{ }
	{}

\swapnumbers

{\theoremstyle{plainsl}
\newtheorem{theorem}{Theorem}[section]
\newtheorem{lemma}[theorem]{Lemma}
\newtheorem{corollary}[theorem]{Corollary}}
{\theoremstyle{remark}
}

\renewcommand\proof{\noindent\textsl{Proof. }}
\newcommand\sqr[2]{{\vbox{\hrule height.#2pt
    \hbox{\vrule width.#2pt height#1pt \kern#1pt
        \vrule width.#2pt}\hrule height.#2pt}}}
\renewcommand{\sqr}{$\blacksquare$}

\renewcommand\qed{%
	\ifmmode\eqno\sqr
	\else\nolinebreak\ \hfill\sqr\medbreak\fi}

\DeclareMathOperator{\rk}{rk}
\DeclareMathOperator{\tr}{tr}
\DeclareMathOperator{\col}{col}

\DeclareMathOperator{\cM}{\mathcal{M}}

\newcommand\one{{\bf1}}

\newcommand\grp[1]{\langle #1\rangle}

\newcommand\sym[1]{\opk{Sym}(#1)}

\usepackage{blkarray}

\title{Quantum Walks on Embeddings}
\author{Hanmeng Zhan \thanks{Department of Combinatorics and Optimization, University of Waterloo; currently a postdoctoral fellow at the Centre de Recherches Math\'ematiques, Universit\'e de Montr\'eal. \texttt{zhanhanm@crm.umontreal.ca}}}
\begin{document}
\maketitle

\begin{abstract}
We introduce a new type of discrete quantum walks, called vertex-face walks, based on orientable embeddings. We first establish a spectral correspondence between the transition matrix $U$ and the vertex-face incidence structure. Using the incidence graph, we derive a formula for the principal logarithm of $U^2$, and find conditions for its underlying digraph to be an oriented graph. In particular, we show this happens if the vertex-face incidence structure forms a partial geometric design. We also explore properties of vertex-face walks on the covers of a graph. Finally, we study a non-classical behavior of vertex-face walks.

\end{abstract}



\section{Introduction \label{Sec_intro}}
Quantum walks were shown to be universal for quantum computation \cite{Childs2009a,Lovett2010,Underwood2010}. The first model of discrete quantum walks was formally introduced by Aharonov et al \cite{Aharonov2000}. Since then, different models have been studied \cite{Kendon2011,Szegedy2004,Portugal2015} and compared \cite{Wong2016,Konno2018}. 

In this paper, we construct a new discrete quantum walk from an orientable embedding of a graph. Roughly speaking, the walk is defined by two partitions of the arcs: one based on the faces, and one on the vertices. To illustrate the idea, we take the planar embedding of $K_4$ as an example. As shown in Figure \ref{planarK4}, since the surface is orientable, we can choose a consistent orientation of the face boundaries. This partitions the arcs of $K_4$ into four groups $\{f_0,f_1,f_2,f_3\}$, called the \textsl{facial walks}. Meanwhile, the arcs can be partitioned into another four groups, each having the same tail. We represent these two partitions by the incidence matrices in Equation \eqref{incmtxs}.

\begin{figure}[h]
	\centering
\begin{minipage}{0.4\textwidth}
	\begin{tikzpicture}
	[every node/.style={circle, draw}]
	
	\node[] (0) at (0,-0.3) {0};
	\node[] (1) at (0,1.6) {1};
	\node[] (2) at (1.7,-1.6) {2};
	\node[] (3) at (-1.7,-1.6) {3};
	
	\node[draw=none] (f0) at (0.6,0) {$\circlearrowright$};
	\node[draw=none] (f3) at (-0.6,0) {$\circlearrowright$};
	\node[draw=none] (f2) at (0,-1) {$\circlearrowright$};
	\node[draw=none] (f3) at (1.5,0) {$\circlearrowright$};
	
	\foreach \a/\b in {0/1,1/2, 2/3,0/3,1/3,0/2}
	\draw[thick] (\a) to (\b);
	\end{tikzpicture}
\end{minipage}%
\begin{minipage}{0.3\textwidth}
	Facial walks:
	\begin{align*}
	f_0&=\{(0,1),(1,2),(2,0)\}\\
	f_1&=\{(1,3),(3,2),(2,1)\}\\
	f_2&=\{(0,2),(2,3),(3,0)\}\\
	f_3&=\{(0,3),(3,1),(1,0)\}
	\end{align*}
\end{minipage}
\caption{A Planar Embedding of $K_4$}
\label{planarK4}
\end{figure}

\begin{equation}
M = \begin{blockarray}{rrrrr}
& f_0 & f_1 & f_2 & f_3\\
\begin{block}{r(rrrr)}
(0,1) &1 & 0 & 0 & 0 \\
(0,2) & 0 & 0 & 1 & 0 \\
(0,3) &0 & 0 & 0 & 1 \\
(1,0) & 0 & 0 & 0 & 1 \\
(1,2) &1 & 0 & 0 & 0 \\
(1,3) &0 & 1 & 0 & 0 \\
(2,0) & 1 & 0 & 0 & 0 \\
(2,1) &0 & 1 & 0 & 0 \\
(2,3) & 0 & 0 & 1 & 0 \\
(3,0) &0 & 0 & 1 & 0 \\
(3,1) &0 & 0 & 0 & 1 \\
(3,2) &0 & 1 & 0 & 0\\
\end{block}
\end{blockarray}
\qquad \qquad
N = \begin{blockarray}{rrrrr}
& 0 & 1 & 2 & 3\\
\begin{block}{r(rrrr)}
(0,1) &1 & 0 & 0 & 0 \\
(0,2) &1 & 0 & 0 & 0 \\
(0,3) &1 & 0 & 0 & 0 \\
(1,0) &0 & 1 & 0 & 0 \\
(1,2) &0 & 1 & 0 & 0 \\
(1,3) &0 & 1 & 0 & 0 \\
(2,0) &0 & 0 & 1 & 0 \\
(2,1) &0 & 0 & 1 & 0 \\
(2,3) &0 & 0 & 1 & 0 \\
(3,0) &0 & 0 & 0 & 1 \\
(3,1) &0 & 0 & 0 & 1 \\
(3,2) &0 & 0 & 0 & 1\\
\end{block}
\end{blockarray}\label{incmtxs}
\end{equation}

If $\widehat{M}$ is the matrix obtained from $M$ by scaling each column to a unit vector, then $\widehat{M}\widehat{M}^T$ is the orthogonal projection onto $\col(M)$, and so 
\[2\widehat{M}\widehat{M}^T-I\]
is the reflection about $\col(M)$. Similarly, if $\widehat{N}$ is the normalized arc-tail incidence matrix, then
\[2\widehat{N}\widehat{N}^T-I\]
is the reflection about $\col(N)$. Now
\[U:= (2\widehat{M}\widehat{M}^T-I)(2\widehat{N}\widehat{N}^T-I)\]
is a unitary matrix, which serves as the transition matrix of our discrete quantum walk.

We can easily generalize the above construction to any orientable embedding of a graph $X$; such a walk will be called a \textsl{vertex-face walk}. While this model has never been studied, there are search algorithms that effectively use the vertex-face walk of a toroidal embedding of $C_n\square C_n$ \cite{Petel2004,Falk2013,Ambainis2013}. Section \ref{Sec_lit1} will discuss this connection in more detail.

\begin{figure}[h]
\centering
\includegraphics[width=8cm]{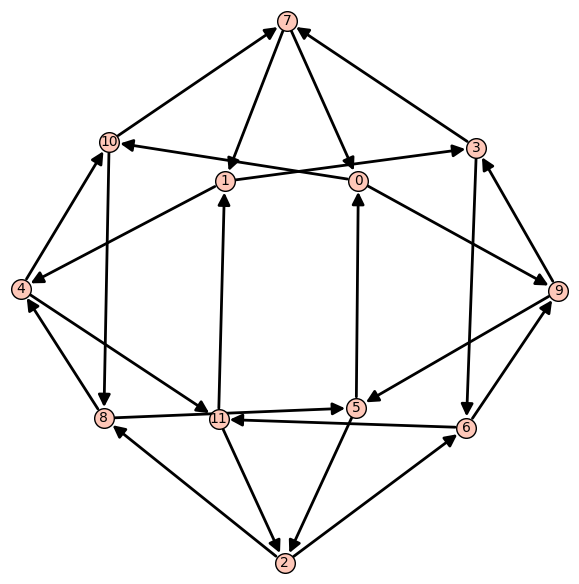}
\caption{$H$-digraph of the planar embedding of $K_4$}
\label{Fig_K4}
\end{figure}

The transition matrix $U$ of the previous example has an interesting expression: $U=\exp(tS)$ where $t\in \re$ and $S$ is the skew-adjacency matrix of an oriented graph, as shown in Figure \ref{Fig_K4}. Thus, $U$ can be seen as the transition matrix of a continuous quantum walk, evaluated at time $t$. We would like to characterize vertex-face walks that are connected to continuous quantum walks in this way. 

Our approach is spectral. We first list some basic properties of the transition matrix of a vertex-face walk. Then we establish a spectral correspondence between the transition matrix and the vertex-face incidence matrix. Using the incidence graph, we derive a formula for the principal logarithm of $U^2$. We then explore necessary conditions and sufficient conditions for the underlying digraph of this logarithm to be an oriented graph, and find interesting connections to partial geometric designs. We also investigate properties of vertex-face walks on the covers of a graph. Finally, we note that some vertex-face walks are reluctant to leave its initial state, which is in sharp contrast to classical random walks. 

\section{Model \label{Sec_model}}
An embedding is \textsl{circular} if every face is bounded by a cycle. In this section, we generalize the example in Section \ref{Sec_intro} to a vertex-face walk on any circular orientable  embedding. 

Let $X$ be a graph, and $\cM$ an embedding of $X$ on some orientable surface. Consider a consistent orientation of the faces, that is, for each edge $e$ shared by two faces $f$ and $h$, the direction $e$ receives in $f$ is opposite to the direction it receives in $h$. Given such an orientation, every arc belongs to exactly one face; let $M$ be the associated arc-face incidence matrix. We also partition the arcs according to their tails, and let $N$ be the associated arc-tail incidence matrix. Denote the normalized versions of $M$ and $N$ by $\widehat{M}$ and $\widehat{N}$. The unitary matrix 
\[U:= (2\widehat{M}\widehat{M}^T-I)(2\widehat{N}\widehat{N}^T-I)\]
is the transition matrix of a vertex-face walk for $\cM$.

Although $U$ depends on the consistent orientation, there are only two choices---reversing all the arcs in the facial walks of one orientation produces the other. Let $R$ be the permutation matrix that swaps arc $(u,v)$ with arc $(v,u)$. If 
\[(2\widehat{M}\widehat{M}^T-I)(2\widehat{N}\widehat{N}^T-I)\]
is the transition matrix relative to the ``clockwise" orientation, then 
\[R(2\widehat{M}\widehat{M}^T-I)R(2\widehat{N}\widehat{N}^T-I)\]
is the transition matrix relative to the ``counterclockwise" orientation. In this paper, we will not specify the orientation when proving properties of $U$, as our results are independent of the choice.

The following observation on duality is immediate.

\begin{lemma}
	If $U$ is the vertex-face transition matrix for $\cM$, then $U^T$ is the vertex-face transition matrix for the dual embedding of $\cM$.
	\qed
\end{lemma}

Define two matrices
\[P:=\widehat{M}\widehat{M}^T,\quad Q:=\widehat{N}\widehat{N}^T.\]
Note that $P$ is the projection onto vectors that are constant on each facial walk, and $Q$ is the projection onto vectors that are constant on arcs leaving each vertex. Let $\one$ denote the all-ones vector. The projections $P$ and $Q$ satisfy the following properties.

\begin{lemma}\label{Lem_props}
For any arc $(u,v)$, let $f_{uv}$ denote the facial walk using $(u,v)$. For any two faces $f$ and $h$, let $f\cap h$ denote the set of vertices used by both $f$ and $h$.
\begin{enumerate}[(i)]
	\item The projections $P$ and $Q$ are doubly stochastic, and so
	\[U\one = U^T \one =\one.\]
	\item For two arcs $(a,b)$ and $(u,v)$, 
	\[P_{(a,b),(u,v)} = \begin{cases}
	\dfrac{1}{\deg(f_{uv})},\quad &\text{if $f_{ab}=f_{uv}$}.\\
	0, \quad &\text{otherwise}.
	\end{cases},\]
	and
	\[Q_{(a,b),(u,v)} = \begin{cases}
	\dfrac{1}{\deg(u)}, \quad & \text{if $a=u$}.\\
	0, \quad & \text{otherwise}.
	\end{cases}\]
	\item For two arcs $(a,b)$ and $(u,v)$, 
	\[(PQ)_{(a,b),(u,v)} = (QP)_{(u,v),(a,b)}=
	\begin{cases}
	\dfrac{1}{\deg(u)\deg(f_{ab})},\quad &\text{ if } u \in f_{ab}.\\
	0, &\text{ otherwise}.
	\end{cases}\]

	\item For two faces $f$ and $h$,
	\[(\widehat{M}^TQ\widehat{M})_{f,h}= \frac{1}{\sqrt{\deg(f)\deg(h)}} \sum_{u \in f\cap h} \frac{1}{\deg(u)}.\]
	For two vertice $u$ and $v$,
	\[(\widehat{N}^TP\widehat{N})_{u,v}=\frac{1}{\sqrt{\deg(u)\deg(v)}}\sum_{f: u, v\in f}\frac{1}{\deg(f)}. \]
\end{enumerate}
\end{lemma}
\proof
We prove the first parts of (iii) and (iv); the other statements follow similarly. Note that $M_{(a,b),f}\ne 0$ if and only if $f=f_{ab}$,  $(M^TN)_{f,w}\ne 0$ if and only if $w$ is contained in $f$, and $N^T_{w, (u,v)}\ne 0$ if and only if $w=u$. Hence $(PQ)_{(a,b),(u,v)}\ne 0$ if and only if $u\in f_{ab}$. This proves (iii). For (iv), simply note that $(M^TNN^TM)_{f,h}$ counts the vertices that appear in both faces $f$ and $h$.
\qed

The above lemma allows us to write out the entries of $U$ explicitly. Moreover, if either $X$ or its dual graph is regular, we have a simple expression for $\tr(U)$. 

\begin{lemma}\label{trace}
Suppose the circular orientable embedding of $X$ has $n$ vertices, $\ell$ edges and $s$ faces. If either $X$ or its dual graph is regular, then
\[\tr(U) = 2 \left( \frac{ns}{\ell} - (n+s-\ell)\right).\]
\end{lemma}
\proof
We have
\[U = (2P-I)(2Q-I).\]
where $P$ and $Q$ are projections. From (iii) in Lemma \ref{Lem_props} we see that
\begin{align*}
\tr(PQ)
&=\sum_{(u,v)} \frac{1}{\deg(u)} \frac{1}{\deg(f_{uv})}\\
&=\sum_f \frac{1}{\deg(f)} \sum_{u \in f} \frac{1}{\deg(u)}.
\end{align*}
If $X$ is $d$-regular, then
\[\tr(PQ) = \frac{s}{d} = \frac{ns}{2\ell}.\]
Hence
\begin{align*}
\tr(U)
&=4\tr(PQ) - 2\tr(P)-2\tr(Q)-\tr(I)\\
&=2\frac{ns}{\ell} - 2(\rk(P) + \rk(Q) - 2\ell)\\
&=2\left(\frac{ns}{\ell} - (n+s-\ell)\right).
\end{align*}
The case where the dual graph is regular follows from $\tr(U)=\tr(U^T)$.
\qed

A quantum walk is called \textsl{reducible} if $U$ is permutation similar to some block-diagonal matrix, and \textsl{irreducible} otherwise. The following result shows that for a connected graph, any vertex-face walk is irreducible, regardless of the embedding.

\begin{lemma}
Let $\cM$ be an orientable embedding of a connected graph $X$. Let $\pi_1$ and $\pi_2$ be the arc-face partition and the arc-tail partition of $\cM$. Then $\pi_1\wedge \pi_2$ is the discrete partition, an $\pi_1\vee \pi_2$ is the trivial partition.
\end{lemma}
\proof
First of all, since every face is bounded by a cycle, no two arcs sharing the tail are contained in the same facial walk, so $\pi_1\wedge \pi_2$ is the discrete partition. Next, since $X$ is connected, between any two vertices $v_0$ and $v_k$ there is a path, say
\[\seq{v}{0}{1}{k}.\]
Consider the first two arcs $(v_0, v_1)$ and $(v_1, v_2)$. If they belong to the same facial walk, then they are in the same class of $\pi_1\vee \pi_2$. Otherwise, there is an arc $(v_1, w_1)$ that is in the same facial walk as $(v_0, v_1)$. Thus, all outgoing arcs of $v_1$, including $(v_1, v_2)$, are in the same class of $\pi_1\vee \pi_2$ as $(v_0, v_1)$. Proceeding in this fashion, we see that all arcs in the path belong to the same class of $\pi_1\vee \pi_2$.
\qed

\section{Spectral Decomposition \label{Sec_sp}}
We now derive the spectral decomposition of a vertex-face transition matrix $U$. Note that $U$ is a product of two reflections, as
\[(2P-I)^2 = (2Q-I)^2=I.\]
The spectrum of a two-reflection unitary matrix was first studied by Szegedy \cite{Szegedy2004}. Later, Godsil developed some machinery, in his unpublished notes  \cite{Godsil2015a}, towards finding the spectral decomposition of any matrix lying in the algebra generated by two reflections. For details of his treatment, see \cite[Section 2.3]{Zhan2018}; our characterization of the eigenspaces of $U$ can be seen as a consequence of those results.

There are three classes of the eigenvalues of $U$: the real eigenvalue $1$, the real eigenvalue $-1$, and other complex eigenvalues on the unit circle. For each class, we characterize the associated eigenspace using matrices $\widehat{M}$ and $\widehat{N}$. To start, suppose the embedding $\cM$ has $n$ vertices, $\ell$ edges and $s$ faces. Let $g$ be the genus of the orientable surface.

\begin{theorem}\label{Thm_vf_1es}
The $1$-eigenspace of $U$ is
\[(\col(M)\cap \col(N)) \oplus (\ker(M^T)\cap \ker(N^T))\]
with dimension $\ell+2g$. Moreover, the first subspace is simply
\[\col(M)\cap \col(N) = \mathrm{span}\{\one\}.\]
\end{theorem}
\proof
We first prove the structure of the $1$-eigenspace. Clearly, every vector in 
\[(\col(M)\cap \col(N)) \oplus (\ker(M^T)\cap \ker(N^T))\]
is an eigenvector for $U$ with eigenvalue $1$. Now suppose $Uz=z$ for some vector $z$. Then 
\[(2Q-I)z = (2P-I)z.\]
Thus $Pz=Qz$ and $(I-P)z=(I-Q)z$. From the decomposition
\[z = Pz + (I-P)z,\]
we see that $z$ lies in 
\[(\col(M)\cap \col(N)) \oplus (\ker(M^T)\cap \ker(N^T)).\]

Next, note that any vector lying in $\col(M)\cap \col(N)$ must be constant over the arcs leaving each vertex, as well as constant on the arcs used by each face. Since $X$ is connected, this vector is constant everywhere. 

Finally, for the multiplicity, we have
\begin{align*}
\dim(\ker(P)\cap \ker(Q))
&= \dim\left( \ker\pmat{P\\Q}\right)\\
&=2\ell -\rk\pmat{P & Q}\\
&=2\ell -\dim(\col(P)+\col(Q))\\
&=2\ell - (\rk(P)+\rk(Q)-\dim(\col(P)\cap\col(Q)))\\
&=2\ell - n - s + 2\\
&=\ell + 2g. \tag*{\sqr}
\end{align*}

To characterize the remaining eigenspaces for $U$, we introduce a few more incidence matrices. A vertex is incident to a face if it is incident to an edge that is contained in the face. Let $B$, $C$ and $D$ be the vertex-edge incidence matrix, the vertex-face incidence matrix, and the face-edge incidence matrix, respectively. Since every face is bounded by a cycle, we have the following two expressions for $C$. 

\begin{lemma}
The vertex-face incidence matrix $C$ satisfies
\[C = \frac{1}{2}BD^T=N^TM.\tag*{\sqr}\]
\end{lemma}

We also define
\[\widehat{C}:=\widehat{N}^T\widehat{M},\]
and call it the \textsl{normalized vertex-face incidence matrix}. All other eigenspaces for $U$ are determined by $\widehat{C}$.

\begin{theorem}\label{Thm_vf_-1es}
The $(-1)$-eigenspace for $U$ is
\[\widehat{M}\ker(\widehat{C})\oplus \widehat{N}^T\ker(\widehat{C}^T)\]
with dimension 
\[n+s - 2\rk(C). \]
\end{theorem}
\proof
One can apply the same argument as the proof of Theorem \ref{Thm_vf_1es}, by replacing $Q$ with $I-Q$, to show that the $(-1)$-eigenspace for $U$ is
\[(\col(P)\cap \ker(Q)) \oplus (\ker(P)\cap \col(Q)).\]
We prove that $\col(P)\cap \ker(Q)$ is isomorphic to $\widehat{M}\ker(\widehat{C})$; the other isomorphism is similar.

If 
\[\widehat{C}y =0,\]
then
\[Q\widehat{M}y = \widehat{N}\widehat{C}y=0.\]
Hence 
\[\widehat{M}y \in \col(P) \cap \ker(Q).\]
Moreover, since $\widehat{M}$ has full column rank, this map is injective. On the other hand, for any $x\in\col(P) \cap \ker(Q)$, there is some $y$ such that
\[x = \widehat{M}y\]
and
\[0 = Qx = \widehat{N}\widehat{C}y,\]
which implies that 
\[y \in \ker(\widehat{C}).\tag*{\sqr}\]

By interlacing, the eigenvalues of $\widehat{C}\widehat{C}^T$ lie in $[0,1]$. We will show that the eigenspaces for $\widehat{C}\widehat{C}^T$ with eigenvalues in $(0,1)$ give rise to eigenspaces for $U$ with non-real eigenvalues.

\begin{theorem}\label{Thm_vf_nonreales}
The multiplicities of the non-real eigenvalues of $U$ sum to $2\rk(C)-2$. Let $\mu\in(0,1)$ be an eigenvalue of $\widehat{C}\widehat{C}^T$. Choose $\theta$ with 
\[\cos(\theta)=2\mu-1.\]
The map
\[y\mapsto (\cos(\theta)+1)\widehat{N}y - (e^{i\theta}+1)\widehat{M}\widehat{C}^Ty\]
is an isomorphism from the $\mu$-eigenspace of $\widehat{C}\widehat{C}^T$ to the $e^{i\theta}$-eigenspace of $U$, and the map
\[y\mapsto(\cos(\theta)+1)\widehat{N}y - (e^{-i\theta}+1) \widehat{M}\widehat{C}^Ty\]
is an isomorphism from the $\mu$-eigenspace of $\widehat{C}\widehat{C}^T$ to the $e^{-i\theta}$-eigenspace of $U$.
\end{theorem}
\proof
Let $y$ be an eigenvector for $\widehat{C}\widehat{C}^T$ with eigenvalue $\mu$. Set $z = \widehat{N}y$. 
\begin{enumerate}[(i)]
\item Suppose $\mu=1$. Then
\[z^*(I-P)z=0,\]
and it follows from the positive-definiteness of $I-P$ that 
\[z\in \col(P)\cap \col(Q).\]
\item Suppose $\mu=0$. By a similar argument,  \[z\in\ker(P)\cap \col(Q).\]
\item Suppose $0<\mu<1$. Then the subspace spanned by $\{z, Pz\}$ is $U$-invariant:
\[U\pmat{z & Pz} = \pmat{z & Pz} \pmat{-1 & -2\mu \\ 2 & 4\mu-1}.\]
To find linear combinations of $z$ and $Pz$ that are eigenvectors of $U$, we diagonalize the matrix
\[\pmat{-1 & -2\mu \\ 2 & 4\mu-1}.\]
It has two eigenvalues: $e^{i\theta}$ with eigenvector 
\[\pmat{-\cos(\theta)-1\\e^{i\theta}+1},\]
and $e^{-i\theta}$ with eigenvector
\[\pmat{-\cos(\theta)-1\\e^{-i\theta}+1}.\]
Since $0<\mu<1$, these two eigenvalues are distinct, and 
\[\frac{\cos(\theta)+1}{e^{\pm i\theta}+1}I - P\]
is invertible, so 
\[(\cos(\theta)+1)z - (e^{\pm i\theta}+1)Pz\]
is indeed an eigenvector for $U$ with eigenvalue $e^{\pm i\theta}$.
\end{enumerate}
Thus, the eigenvectors for $\widehat{C}\widehat{C}^T$ with eigenvalues in $(0,1)$ provide
\[2(\rk(C) - \dim(\col(P) \cap \col(Q))) = 2\rk(C)-2\]
orthogonal eigenvectors for $U$, which span the orthogonal complement of the $(\pm1)$-eigenspaces.
\qed

After normalization, we obtain an explicit formula for the eigenprojection of each non-real eigenvalue of $U$.

\begin{corollary}\label{Cor_vf_nonrealeproj}
	Let $\mu\in(0,1)$ be an eigenvalue of $\widehat{C}\widehat{C}^T$. Choose $\theta$ such that $\cos(\theta)=2\mu-1$. Let $E_{\mu}$ be the orthogonal projection onto the $\mu$-eigenspace of $\widehat{C}\widehat{C}^T$. Set
	\[W:=\widehat{N}E_{\mu}\widehat{N}^T.\]
	Then the $e^{i\theta}$-eigenprojection of $U$ is
	\[\frac{1}{\sin^2(\theta)}\left((\cos(\theta)+1)W-(e^{i\theta}+1)PW-(e^{-i\theta}+1)WP+2PWP\right),\]
	and the $e^{-i\theta}$-eigenprojection of $U$ is
	\[\frac{1}{\sin^2(\theta)}\left((\cos(\theta)+1)W-(e^{-i\theta}+1)PW-(e^{i\theta}+1)WP+2PWP\right). \tag*{\sqr}\]
\end{corollary}

\section{Hamiltonian \label{Sec_H}}
A \textsl{Hamiltonian} of a unitary matrix $V$ is a Hermitian matrix $H$ such that $V=\exp(iH)$. Given the spectral decomposition
\[V = \sum_r \alpha_r F_r,\]
any Hamiltonian $H$ can be written as
\[H = -i \sum_r \log(\alpha_r) F_r,\]
for some value of $\log(\alpha_r)$. If in addition, for each $r$, the angle satisfies 
\[-\pi < -i\log(\alpha_r) \le \pi,\]
then $H$ is called the \textsl{principal Hamiltonian of $V$}. In this sense, every unitary matrix can be viewed as the transition matrix of a continuous quantum walk on the underlying digraph of its principal Hamiltonian.

We study the principal Hamiltonian of $U^2$, where $U$ is the transition matrix of a vertex-face walk. The spectral machinery we developed in the last section reveals a close connection between $H$ and the bipartite vertex-face incidence graph. For simplicity, we will focus on circular orientable embeddings where both $X$ and the dual graph are regular; if each vertex has $d$ neighbors and each face uses $k$ vertices, we say the embedding \textsl{has type $(k,d)$}. 

\begin{theorem}\label{Thm_H}
Let $\cM$ be a circular orientable embedding of type $(k,d)$. Let $M$ and $N$ be the arc-face incidence matrix and the arc-tail incidence matrix, respectively. Let $U$ be the transition matrix of a vertex-face walk for $\cM$. Let $A$ be the adjacency matrix of the vertex-face incidence graph, with spectral decomposition
\[A = \sum_{\lambda} \lambda G_{\lambda}.\]
Then each eigenvalue $\lambda$ of $A$ gives rise to some eigenvalue $e^{\pm i\theta}$ of $U$ in the following way:
\[\frac{\abs{\lambda}}{\sqrt{dk}}=\cos\left(\frac{\theta}{2}\right).\]
Moreover, $U^2=\exp(iH)$, where 
\[H = 4 \pmat{N & iM} \left(\sum_{\lambda\notin\{0, \pm\sqrt{dk}\}} \frac{\lambda\arccos(\abs{\lambda}/\sqrt{dk})}{\abs{\lambda}\sqrt{dk-\lambda^2}} G_r\right)\pmat{N^T\\-iM^T}.\]
\end{theorem}
\proof
As before, let $C$ be the vertex-face incidence matrix. Then
\[A = \pmat{ 0 & C\\ C^T & 0}.\]
Note that the eigenvalues of $A$ are symmetric about zero, and bounded in absolute value by $\sqrt{dk}$. By Theorems \ref{Thm_vf_1es}, \ref{Thm_vf_-1es} and \ref{Thm_vf_nonreales}, each eigenvalue $\lambda$ of $A$ determines the real part of some eigenvalue $e^{\pm i \theta}$ of $U$:
\[\cos(\theta)=\frac{2\lambda^2}{dk}-1,\]
that is,
\[\frac{\abs{\lambda}}{\sqrt{dk}}=\cos\left(\frac{\theta}{2}\right).\]

Moreover, for $\lambda\ne 0$, the $\lambda$-eigenprojection for $A$ is 
\[G_{\lambda}=\frac{1}{2}\pmat{E_{\lambda^2} &  \frac{1}{\lambda}E_{\lambda^2}C\\ \frac{1}{\lambda}C^T E_{\lambda^2} & \frac{1}{\lambda^2}C^TE_{\lambda^2} C},\]
where $E_{\lambda^2}$ is the $\lambda^2$-eigenprojection for $CC^T$. Since $U$ has real entries, its spectrum is closed under complex conjugation. Hence by Corollary \ref{Cor_vf_nonrealeproj},
\[H = -4i\sum_{0<\lambda<\sqrt{dk}}\frac{\arccos(\lambda/\sqrt{dk})}{\lambda \sqrt{dk-\lambda^2}}(NE_{\lambda^2}CM^T-MC^TE_{\lambda^2}N^T).\]
Combining this with the fact that
\[G_{\lambda}-G_{-\lambda} = \pmat{0 & \frac{1}{\lambda}E_{\lambda^2} C\\ \frac{1}{\lambda}C^T E_{\lambda^2} & 0}\]
yields the formula for $H$.
\qed

The above result implies that $H$ is the block sum of
\[\pmat{N & 0 \\ 0 & iM} \phi(A) \pmat{N^T & 0\\0&-iM^T}\]
for some odd polynomial $\phi$.

\section{$H$-digraph \label{Sec_Hdigraph}}
Given an embedding $\cM$ and the vertex-face transition matrix $U$, we will refer to the underlying digraph of the principal Hamiltonian $H$ of $U^2$ as the \textsl{$H$-digraph} of $\cM$. By Theorem \ref{Thm_H}, $iH$ is skew-symmetric, so the $H$-digraph is a weighted oriented graph. 

An embedding $\cM$ is \textsl{orientably-regular} if its orientation-preserving automorphism group acts regularly on the arcs. Using the decomposition 
\[C = \frac{1}{2}BD^T=N^TM,\]
where $B$, $C$, $D$, $M$, $N$ are the vertex-edge,  vertex-face, face-edge, arc-face, and arc-tail incidence matrices, we obtain the following.

\begin{theorem}
Let $\cM$ be a circular orientable embedding of type $(k,d)$. If $\cM$ is orientably-regular, then the vertex-face incidence graph is edge-transitive, and the $H$-digraph is vertex-transitive.
\qed
\end{theorem}

In general, we should expect the $H$-digraph to be dense with many different weights; however, there are cases where it is sparse and unweighted (up to scaling). An example was given in Figure \ref{Fig_K4}. In this section, we study circular orientable embeddings of type $(k,d)$ whose $H$-digraphs are oriented graphs. We first give a necessary condition on the eigenvalues of the vertex-face incidence graph. This is a direct consequence of Godsil's observation on real state transfer \cite{Godsil2017b}, which we summarize below.

\begin{theorem}\label{Thm_ratio}
Let $H$ be a Hermitian matrix with algebraic entries. Suppose for some real number $t$ the entries of $\exp(itH)$ are algebraic. Then the ratio of any two non-zero eigenvalues of $H$ are rational. Moreover, if $iH$ has integer entries, then there is a square-free integer $\Delta$ such that all eigenvalue of $H$ are in $\ints[\sqrt{\Delta}]$.
\qed
\end{theorem}

The ratio condition in Theorem \ref{Thm_ratio} is particularly useful in characterizing state transfer in continuous quantum walks on graphs (see for example \cite{Godsil2011a}) and oriented graphs (\cite{Godsil2017b,Lato2019}). Here, we present its application to discrete quantum walks.

\begin{theorem}
Let $U$ be a vertex-face transition matrix for a circular orientable embedding of type $(k,d)$. Let $A$ be the adjacency matrix of the vertex-face incidence graph.  Suppose $U^2=\exp(tS)$ for some real number $t$ and skew-adjacency matrix $S$. Then the following hold.
\begin{enumerate}[(i)]
\item There is a square-free integer $\Delta$ such that all eigenvalues of $S$ are in $\ints[\sqrt{-\Delta}]$.
\item If $\lambda_r$ and $\lambda_s$ are two eigenvalues of $A$ that are distinct from $\{0,\pm\sqrt{dk}\}$, then
\[\frac{\arccos(\abs{\lambda_r}/\sqrt{dk})}{\arccos(\abs{\lambda_s}/\sqrt{dk})}\in \rats.\]
\end{enumerate}
\end{theorem}
\proof 
(i) follows from Theorem \ref{Thm_H} since $U$ has rational entries and $S$ has integer entries. For (ii), recall that $\lambda_r$ and $\lambda_s$ determine non-real eigenvalues $e^{\pm i\theta_r}$ and $e^{\pm i\theta_s}$ of $U$ by
\[\frac{\abs{\lambda_r}}{\sqrt{dk}}=\cos\left(\frac{\theta_r}{2}\right),\quad \frac{\abs{\lambda_s}}{\sqrt{dk}}=\cos\left(\frac{\theta_s}{2}\right).\]
If $0<\theta_r, \theta_s<\pi$, then $\theta_r/\theta_s$ equals the ratio of two eigenvalues of $S$, which must be rational.
\qed

The above condition is satisfied when $U^2$ has exactly three eigenvalues. In fact, the $H$ digraph of such $U^2$ is guaranteed to be an oriented graph.

\begin{theorem}\label{Thm_oriented}
Let $\cM$ be a circular orientable embedding of type $(k,d)$ of a graph $X$. Let $U$ be a vertex-face transition matrix for $\cM$. Then
\[U^2=\exp(\gamma(U-U^T))\]
for some real number $\gamma$ if and only if the vertex-face incidence graph has four or five distinct eigenvalues. Moreover, 
\[\frac{dk}{4}(U^T-U)\] 
is the skew-adjacency matrix of some oriented graph on the arcs of $X$, and the degree of $(a,b)$ is 
\[dk - \sum_{u \in f_{ab}}\beta(a,u),\]
where $f_{ab}$ denotes the unique face using the arc $(a,b)$, and $\beta(a, u)$ denotes the number of faces containing both $a$ and $u$. 
\end{theorem}
\proof
Let 
\[U = \sum_r \alpha_r F_r\]
be the spectral decomposition of $U$. Then
\[U^2=\exp(\gamma(U-U^T))\]
holds if and only if 
\[\sum_r \alpha_r^2 F_r =\sum_r e^{\gamma(\alpha_r-\alpha_r^{-1})} F_r,\]
that is, for each non-real eigenvalue $\alpha_r = e^{i\theta_r}$ of $U$,
\[e^{2i\theta} = e^{2\gamma \sin(\theta)}.\]
Since $\sin(x)/x$ is monotone when $0<x<\pi$, the above holds for some $\gamma$ if and only if $U^2$ has exactly three eigenvalues, or equivalently, the vertex-face incidence graph has exactly four or five eigenvalues.

Let $M$ and $N$ be the arc-face and arc-tail incidence matrices, respectively. Recall that
\[U = \left(\frac{2}{k}MM^T-I\right)\left(\frac{2}{d}NN^T-I\right).\]
Thus
\[U-U^T= \frac{4}{dk} (MM^TNN^T-NN^TMM^T).\]
Let $S= MM^TNN^T-NN^TMM^T$. By Lemma \ref{Lem_props} (iii), 
\[ S_{(a,b),(u,v)}=
\begin{cases}
1,\quad \text{ if $u\in f_{ab}$ and $a\notin f_{uv}$},\\
-1,\quad \text{ if $a\in f_{uv}$ and $u\notin f_{ab}$},\\
0,\quad \text{ otherwise}.
\end{cases}\]
Therefore $S$ is the skew-adjacency matrix of some oriented graph. Moreover, for each $u\in f_{ab}$, there is a bijection between the neighbors of $v$ such that $a\in f_{uv}$ and the faces containing both $a$ and $u$. Hence the degree of the oriented graph is 
\[\sum_{u\in f_{ab}}(d-\beta(a,u)) = dk - \sum_{u \in f_{ab}} \beta(a,u). \tag*{\sqr}\]

Let $C$ be the vertex-face incidence matrix of the embedding in Theorem \ref{Thm_oriented}. We see that $U^2=\exp(\gamma(U-U^T))$ if and only if $C$ has exactly two non-zero singular values. Combinatorial designs with two non-zero singular eigenvalues were studied by van Dam and Spence \cite{vanDam2004, vanDam2005}. In particular, they showed that a point-$d$-regular and block-$k$-regular design with two non-zero singular values is a \textsl{partial geometric design} with parameters $(d,k,t,c)$, originally introduced by Bose et al \cite{Bose1973}, where for each point-block pair $(p,B)$, the number of incident point-block pairs 
\[\abs{\{(p',B'): p'\ne p, B'\ne B, p'\in B, p\in B'\}}\]
equals $c$ or $t$, depending on whether $p$ is in $B$ or not. Below we include a proof.

\begin{theorem}
Let $C$ be an incidence matrix with $C\one =d\one$ and $C^T\one =k\one$.  Suppose $C$ has exactly two non-zero singular values. Then $C$ is the incidence matrix of a partial geometric design.
\end{theorem}
\proof
Clearly, $dk$ is an eigenvalue of $CC^T$ with eigenprojection $J/n$. Let $\mu$ be the other non-zero eigenvalue of $CC^T$. Let $E_0$ be the projection onto the kernel of $CC^T$. We have
\[CC^T = dk\left(\frac{1}{n}J\right) + \mu \left(I-\frac{1}{n}J-E_0\right)=\mu I+\frac{dk-\mu}{n}J-\mu E_0.\]
Thus
\[CC^TC = \mu C+\frac{k(dk-\mu)}{n} J.\]
Note that $(CC^TC)_{p,B}$ counts all pairs $(p',B')$ with $p'\in B$ and $p\in B'$. Hence the incidence structure is a 
\[\left(d,k,\frac{k(dk-\mu)}n, \frac{k(dk-\mu)}n+\mu+1-d-k\right)\]
partial geometric design.
\qed

We briefly discuss embeddings that realize partial geometric designs. 

If $CC^T$ is invertible, then a trace argument shows that the vertex-face incidence structure is a \textsl{balanced incomplete block design (BIBD)}, or a \textsl{$2$-design}, with parameters $(n,k, d(k-1)/(n-1))$, that is, a point-regular and block-regular design where every two distinct points lie in $d(k-1)/(n-1)$ blocks. 

The study of connections between $2$-designs and graph embeddings dates back to 1897 \cite{Heffter1897}, when Heffter constructed $2$-$(n,3,2)$ designs using certain  triangular embeddings of $K_n$. However, building these triangular embeddings themselves remained a challenging task, until Ringel \cite{Ringel1961}, Gustin \cite{Gustin1963} and Terry et al \cite{Terry1967} provided solutions to all admissible $n$, that is, $n\equiv 0,3,4,7 \pmod{12}$. 

The self-dual circular embeddings of $K_n$, on the other hand, yield a family of $2$-$(n, n-1, n-2)$ designs. While self-dual embeddings of $K_n$ exists if and only if $n\equiv 0,1 \pmod{4}$ \cite{White1984}, the circular ones are only known when $n$ is a prime power. For the constructions, see Biggs \cite{Biggs1971}; we remark that these  are all orientably-regular embeddings.

Embeddings with singular $CC^T$ and $C^TC$ may be related to other designs. A \textsl{two-class partially balanced incomplete block design (PBIBD)} with parameters $(n, k; \lambda_1, \lambda_2)$ is a point-regular, block-regular design whose incidence matrix $C$ satisfies
\[CC^T = dI+\lambda_1 A + \lambda_2(J-I-A),\]
where $A$ is the adjacency matrix of a strongly regular graph. A PBIBD has at most three non-zero singular values; those with two non-zero singular values are precisely partial geometric PBIBDs, and they are usually referred to as \textsl{special PBIBDs} \cite{Bridges1974}. 

Every triangular embedding of a strongly regular graph on $n$ vertices determines a $(n,3;2,0)$-PBIBD. In \cite{Petroelje1971}, Petroelje gave a construction for orientable triangular embeddings of $K_{n,n,n}$, which yield $(3n, 3;2,0)$-special PBIBDs.

\section{Covers \label{Sec_cover}}
In this section, we consider a covering construction that preserves nice properties of vertex-face walks.

An \textsl{arc-function of index} $r$ of $X$ is a map $\phi$ from the arcs of $X$ into $\sym{r}$, such that $\phi(u,v)= \phi(v,u)^{-1}$. The \textsl{fiber} of a vertex $u$ is the set
\[\{(u, i): i=0,1,\cdots, r-1\}.\]
If we replace each vertex of $X$ by its fiber, and join $(u,i)$ to $(v,j)$ whenever $\phi(u,v)(i)=j$, then we obtain a new graph $X^{\phi}$, called the \textsl{$r$-fold cover} of $X$. For example, we can let $\phi$ be the constant arc-function that sends every arc to $(1,2)\in \sym{2}$. Then the double cover $K_4^{\phi}$ is isomorphic to the cube, as shown in Figure \ref{Fig_doublecoverK4}.

\begin{figure}[h]
	\centering
\begin{minipage}[b]{0.4\textwidth}
	\begin{tikzpicture}
	[every node/.style={circle, draw}]
	
	\node[] (0) at (0,-0.3) {0};
	\node[fill=orange!30] (1) at (0,1.6) {1};
	\node[fill=green!30] (2) at (1.7,-1.6) {2};
	\node[fill=red!30] (3) at (-1.7,-1.6) {3};
	
	\foreach \a/\b in {0/1,1/2, 2/3,0/3,1/3,0/2}
	\draw[thick] (\a) to (\b);
	
	\end{tikzpicture}
	\captionof{figure}{$K_4$}
\end{minipage}%
\begin{minipage}[b]{0.45\textwidth}
\begin{tikzpicture}
	[every node/.style={rectangle, draw}]
	
	\node[] (00) at (-1,1) {(0,0)};
	\node[fill=green!30] (21) at (1,1) {(2,1)};
	\node[fill=orange!30] (10) at (1,-1) {(1,0)};
	\node[fill=red!30] (31) at (-1,-1) {(3,1)};
	
	\node[fill=orange!30] (11) at (-2,2) {(1,1)};
	\node[fill=red!30] (30) at (2,2) {(3,0)};
	\node[] (01) at (2,-2) {(0,1)};
	\node[fill=green!30] (20) at (-2,-2) {(2,0)};
	
	\foreach \a/\b in {00/11, 00/21, 00/31, 01/10, 01/20, 01/30,10/21,10/31,11/20,11/30,20/31,21/30}
	\draw[thick] (\a) to (\b);
\end{tikzpicture}
\captionof{figure}{A double cover of $K_4$}
\label{Fig_doublecoverK4}
\end{minipage}
\end{figure}

The above definition tells us how to construct a cover from a base graph. Alternatively, we say a graph $Y$ \textsl{covers} $X$ if there is a homomorphism $\psi$ from $Y$ to $X$, such that for any vertex $y$ of $Y$ and $x=\psi(y)$, the homomorphism $\psi$ restricted to $N_Y(y)$ is a bijection onto $N_X(x)$. The map $\psi$ is called a \textsl{covering map}. If $X$ is connected, then the preimages $\psi^{-1}(x)$ all have the same size; they are precisely the fibers of $X$.

Given an orientable embedding $\cM_X$ of $X$, and a covering map $\psi$ from a connected graph $Y$ to $X$, we define an orientable embedding $\cM_Y$ of $Y$ by specifying its facial walks. Let $W$ be a facial walk of $\cM_X$ starting at vertex $u$. Clearly, the preimage $\psi^{-1}(W)$ consists of walks starting and ending in the fiber $\psi^{-1}(u)$, and each arc of $Y$ appears in at most one of these walks. Then, the facial walks of $\cM_Y$ are exactly the closed walks in the preimages of the facial walks of $\cM_X$. In the previous example, the planar embedding of $K_4$ gives rise to an embedding of the cube on the torus, with $4$ faces each of length $6$.

We will focus on a special type of cover, known as the voltage graphs. A \textsl{voltage graph} of $X$ is an $r$-fold cover $Y=X^{\phi}$, where the image of the arc-function $\phi$ is a subgroup $\Gamma \leq \sym{r}$ of order $r$, and 
\[V(Y) = V(X) \times \Gamma,\quad E(Y) = E(X) \times \Gamma.\]
Voltage graphs correspond to normal covers \cite{Hatcher2002}, and have been extensively studied. We only state one property that voltage graphs satisfy; for more background, see Gross and Tucker \cite{Gross2001}.

\begin{theorem}\label{Thm_lift}
Let $C$ be a $k$-cycle in $X$. Let $Y=X^{\phi}$ be a voltage graph of order $r$. If $\phi(C)$ has order $\ell$, then $C$ lifts to $r/\ell$ cycles in $Y$, each of length $k\ell$.
\qed
\end{theorem}

We call $\cM_Y$ a \textsl{voltage embedding} if $Y$ is a voltage graph of $X$, and $\cM_Y$ is obtained from $\cM_X$ by the above lifting method.

The next result shows that the transition matrix of $\cM_X$ is a block sum of the transition matrix of $\cM_Y$, and consequently, the $H$-digraph of $\cM_X$ is a quotient digraph of $\cM_Y$. To prove it, we need the concept of row and column equitable partitions, which were introduced by Godsil \cite[Ch 12]{Godsil1993}. Let $A$ be a matrix over $\cx$. Let $\sigma$ and $\rho$ be the partition of the columns and rows of $A$, and let $K$ and $L$ be their respective characteristic matrices.  The pair $(\rho,\sigma)$ is \textsl{column equitable} if 
$\col(AK) \sbs \col(L)$, \textsl{row equitable} if $\col(A^*L) \sbs \col(K)$, and \textsl{equitable} if it is both column and row equitable.

\begin{theorem}\label{Thm_vf_quotient}
Let $\cM_X$ be a circular orientable embedding of $X$. Let $Y$ be a voltage graph of $X$, and $\cM_Y$ the associated voltage embedding. Let $\rho$ be the partition of the arcs of $Y$, where each class is the preimage of some arc of $X$. Let $\widehat{L}$ be its normalized incidence matrix of $\rho$. If $U_X$ and $U_Y$ are the vertex-face transition matrices for $\cM_X$ and $\cM_Y$, then 
\[U_X = \widehat{L}^T U_Y \widehat{L}.\]
Consequently, the $H$-digraph of $\cM_X$ is a quotient digraph of the $H$-digraph of $\cM_Y$.
\end{theorem}
\proof
Let $\widehat{M}_X, \widehat{M}_Y, \widehat{N}_X, \widehat{N}_Y$ be the arc-face incidence matrices and arc-tail incidence matrices of the embeddings of $X$ and $Y$, respectively. We have
\[U_Y = (2\widehat{M}_Y \widehat{M}_Y^T - I) (2\widehat{N}_Y \widehat{N}_Y^T-I).\]
Let $\sigma$ be the partition of the vertices of $Y$ into fibers, with normalized incidence matrix $\widehat{K}$. It is not hard to verify that
\[\widehat{N}_Y\widehat{K} = \widehat{L}\widehat{N}_X\]
and
\[\widehat{N}_Y^T \widehat{L} = \widehat{K}\widehat{N}_X^T.\]
Thus $(\rho, \sigma)$ is an equitable partition of $\widehat{N}_Y$. It follows that 
\begin{equation}
\widehat{N}_Y\widehat{K}\widehat{K}^T = \widehat{L}\widehat{L}^T \widehat{N}_Y.
\label{Eqn_eqpartns}
\end{equation}
Since 
\[\widehat{N}_X = \widehat{L}^T\widehat{N}_Y \widehat{K},\]
the projection onto $\col(\widehat{N}_X)$ can be written as
\begin{align*}
\widehat{N}_X\widehat{N}_X^T &= \widehat{L}^T(\widehat{N}_Y\widehat{K}\widehat{K}^T)\widehat{N}_Y^T\widehat{L}\\
&=\widehat{L}^T(\widehat{L}\widehat{L}^T\widehat{N}_Y)\widehat{N}_Y^T\widehat{L}\\
&=\widehat{L}^T\widehat{N}_Y\widehat{N}_Y^T\widehat{L}.
\end{align*}
Applying a similar argument to the preimages of facial walks, we can show that 
\[\widehat{M}_X\widehat{M}_X^T = \widehat{L}^T \widehat{M}_Y\widehat{M}_Y^T \widehat{L}.\]
Thus,
\begin{equation}
U_X = \widehat{L}^T(2\widehat{M}_Y\widehat{M}_Y^T-I)\widehat{L}\widehat{L}^T(2\widehat{N}_Y\widehat{N}_Y^T-I)\widehat{L}.
\label{Eqn_quotient_trans}
\end{equation}
Finally, from Equation \eqref{Eqn_eqpartns} we see that
\[\widehat{L}\widehat{L}^T \widehat{N}_Y\widehat{N}_Y^T=\widehat{N}_Y\widehat{K}\widehat{K}^T\widehat{N}_Y^T,\]
which is a symmetric matrix, so $\widehat{L}\widehat{L}^T$ commutes with $\widehat{N}_Y\widehat{N}_Y^T$. Therefore, Equation \eqref{Eqn_quotient_trans} reduces to 
\[U_X=\widehat{L}^TU_Y\widehat{L}.\tag*{\sqr}\]

Conversely, we may ``lift" nice properties of $\cM_X$ when taking a voltage embedding, using the following simple technique. Let $C$ be the incidence matrix of a design. Construct a new design with incidence matrix $C\otimes \one$ by duplicating the points and preserving the incidence relation. If $C$ is a partial geometric design, then so is $C\otimes \one$. 

\begin{theorem}
Let $\cM_X$ a circular orientable embedding of $X$. Suppose its vertex-face incidence structure is a partial geometric design. Let $Y=X^{\phi}$ be a voltage graph of order $r$, and $M_Y$ the associated voltage embedding. Suppose for each facial cycle $C$ of $\cM_X$, the order of $\phi(C)$ is $r$. Then the vertex-face incidence structure of $\cM_Y$ is also a partial geometric design.
\end{theorem}
\proof
Let $\psi$ be the covering map. By Theorem \ref{Thm_lift}, each facial cycle $f_{ab}$ of $\cM_X$ lifts to a unique facial cycle $\psi^{-1}(f_{ab})$ of $\cM_Y$. Moreover, all arcs in $\psi^{-1}((a,b))$ are contained in $\psi^{-1}(f_{ab})$. 
\qed

This construction yields many new embeddings whose $H$-digraphs are oriented graphs. For example, we have the following family based on the circular self-dual embeddings of $K_n$.

\begin{corollary}
Let $n$ be a power of $2$. Let $\cM_{K_n}$ be a circular self-dual embedding of $K_n$. Let $Y=K_n^{\phi}$ be the double cover of $K_n$ with $\phi$ sending every arc to the involution. Let $\cM_Y$ be the voltage embedding. Then the $H$-digraph of $\cM_Y$ is an oriented graph.
\qed
\end{corollary}

\section{Sedentary Walks}
One counterintuitive phenomenon in quantum walks is that the walker may be reluctant to leave its initial state. This was first observed in continuous quantum walks on $K_n$: for any time $t$, the mixing matrix
\[U_{K_n}(t)\circ \comp{U_{K_n}(t)}\]
converges to $I$ as $n$ goes to infinity. In \cite{Godsil2017a}, Godsil investigated quantum walks on complete graphs, some cones and some strongly regular graphs that enjoy the same property. Following his paper, we say a sequence of discrete quantum walks, determined by transition matrices $\{U_1, U_2, \cdots\}$, is \textsl{sedentary} if for any step $t$, the mixing matrices $U_n^t\circ \comp{U_n^t}$ converges to $I$ as $n$ goes to infinity.

\begin{lemma}\label{Thm_tr}
Let $\cM$ be a circular orientable embedding of type $(d,k)$. Suppose the vertex-face incidence structure is a $2$-design. Then 
\[\tr(U^t)=nd-2(1-\cos(t\theta))(n-1),\]
where 
\[\cos(\theta)=\frac{2(n-k)}{k(n-1)}-1.\]
\end{lemma}
\proof
Let $C$ be the vertex-face incidence matrix. We have
\[CC^T = \frac{d(n-k)}{n-1}I+\frac{d(k-1)}{n-1}J.\]
The eigenvalues of $CC^T$ are $dk$ with multiplicity $1$, and $d(n-k)/(n-1)$ with multiplicity $n-1$. By Theorem \ref{Thm_vf_nonreales}, the non-real eigenvalues of $U$ are $e^{\pm i\theta}$, each with multiplicity $n-1$, where 
\[\cos(\theta)=\frac{2(n-k)}{k(n-1)}-1.\]
Hence $1$ is an eigenvalue of $U$ with multiplicity $nd-2(n-1)$. Therefore,
\[\tr(U^t)=(e^{it\theta}+e^{-it\theta})(n-1)+nd-2(n-1),\]
from which the result follows.
\qed

We found one family of sedentary walks from embeddings we visited before.

\begin{corollary}\label{Thm_sedKn}
	For each prime power $n$, let $U_n$ be the vertex-face transition matrix for a self-dual orientably-regular embedding of $K_n$. The quantum walks determined by
	\[\{U_n: n\text{ is a prime power}\}\]
	form a sedentary family.
\end{corollary}
\proof
Since the embedding is orientably-regular, $U^t$ has constant diagonal. By Theorem \ref{Thm_tr}, each diagonal entry of $U^t \circ U^t$ is
\[\left(1-\frac{2(1-\cos(t\theta))}{n}\right)^2,\]
which converges to $1$ as $n$ goes to infinity. 
\qed

\section{Search\label{Sec_lit1}}
We mention a potential application of the vertex-face walks. First, let us revisit a quantum walk based algorithm due to Patel et al \cite{Petel2004} and Falk \cite{Falk2013}; its performance was proved to match the best known quantum algorithms for searching a marked item on a $2$-dimensional grid \cite{Ambainis2013} .

\begin{figure}[h]
\centering
\includegraphics[width=7cm]{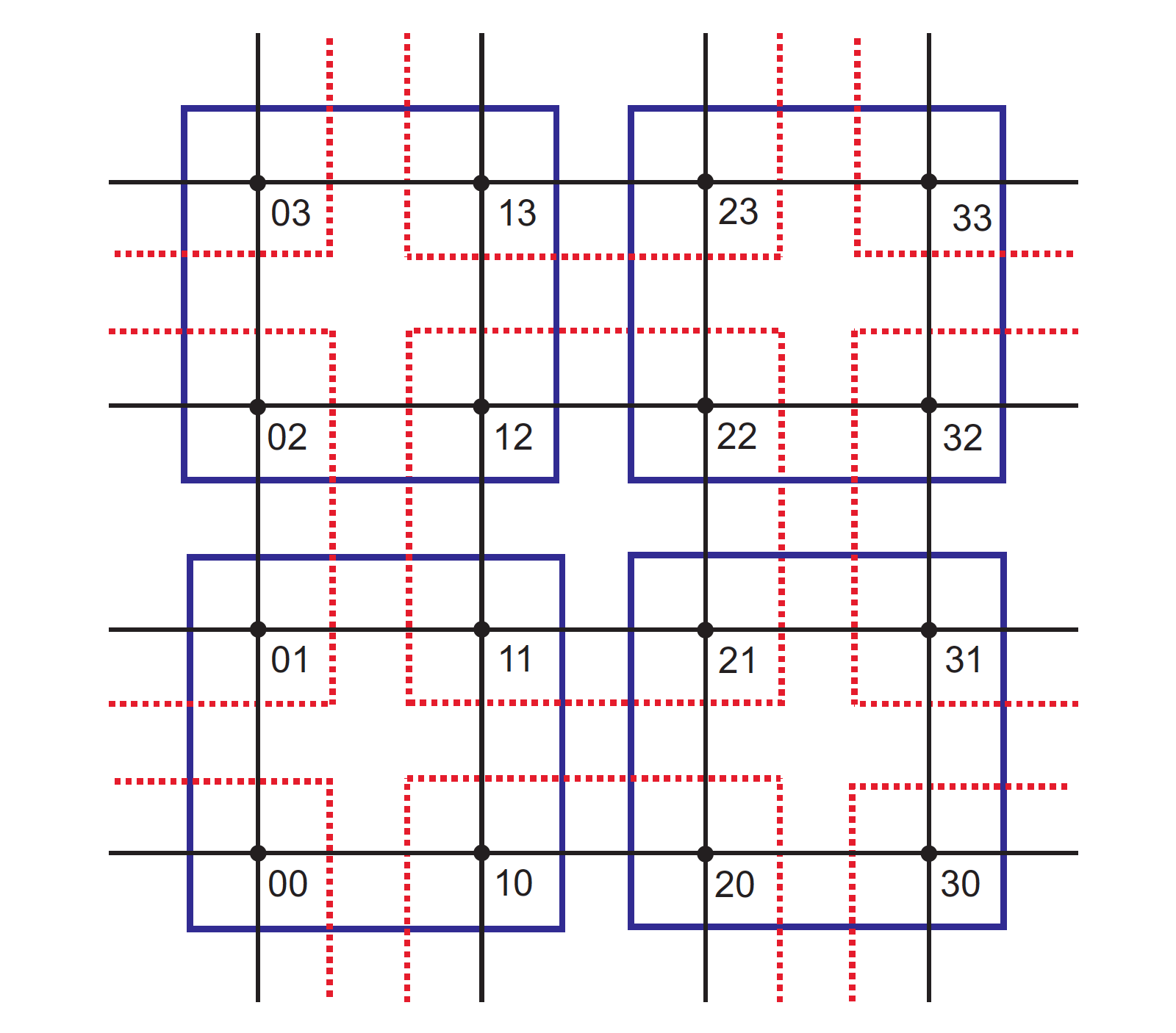}
\caption{Two partitions of the vertices of $C_n \square C_n$}
\label{Fig_grid}
\end{figure}

We will view the $2$-dimensional grid as a Cartesian square of a cycle:
\[X:=C_n \square C_n.\]
Consider two partitions of $V(X)$, illustrated by the blue squares and the red squares in Figure \ref{Fig_grid}. Let $U_o$ and $U_e$ be the reflections about the column spaces of the characteristic matrices of these two partitions, respectively. If we remove the oracle from the search algorithm proposed by \cite{Petel2004,Falk2013}, then it is equivalent to a quantum walk with transition matrix
\[U = U_o U_e.\]

\begin{figure}[h]
\centering
\includegraphics[width=6.5cm]{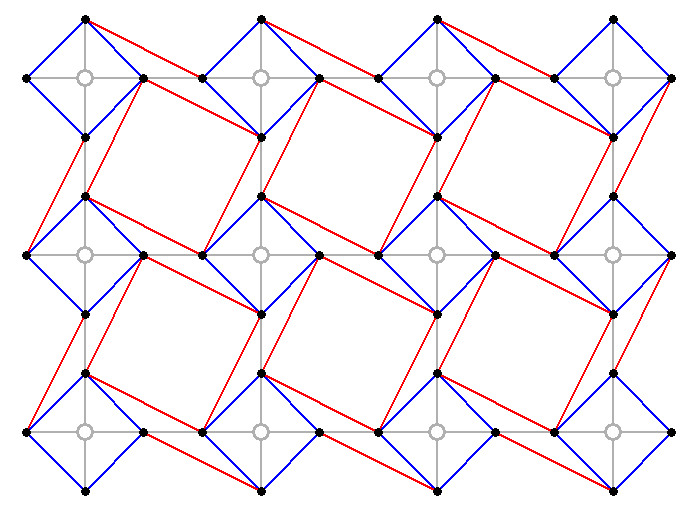}
\caption{Two partitions of the arcs of $C_n \square C_n$}
\label{Fig_arcgrid}
\end{figure}

Notice that Figure \ref{Fig_grid} represents a self-dual embedding of $C_n\square C_n$ on the torus. In fact, it gives rise to a toroidal embedding of another graph $Y$, obtained by truncating the edges of $X$ and joining the new vertices by blue and red edges, as shown in Figure \ref{Fig_arcgrid}.

Clearly, $Y$ is isomorphic to $C_{2n}\square C_{2n}$, and the blue and red squares partition $V(Y)$ the same way they do in Figure \ref{Fig_grid}. Thus based on \cite{Ambainis2013}, we can construct a transition matrix $U$ from these two partitions. On the other hand, we may think of the vertices of $Y$ as arcs of $X$---the one closer to $u$ on edge $\{u,v\}$ is the arc $(u,v)$, with tail $u$. Thus, the blue squares partition the arcs based on their tails, while the red squares partition the arcs based on the faces they lie in. Therefor $U$ is a vertex-face transition matrix for the toroidal embedding of $X$.

In general, given a quantum walk with transition matrix $U$ on the arcs of a graph, we may search for a marked vertex $u$ in the following way. Let $O$ be the matrix that maps $e_a$ to $-e_a$ if $a$ is an outgoing arc of $u$, and fixes $e_a$ otherwise; this is called the \textsl{oracle}. Initialize the system to $\one/\sqrt{m}$, where $m$ is the number of arcs. Apply $OU$ to the initial state $t$ times. The probability of finding $u$ after $t$ steps is given by
\[\sum_{a\text{ has tail } u}\abs{\Big\langle(OU)^t \frac{1}{\sqrt{m}}\one, e_a\Big\rangle}^2.\]

An \textsl{arc-reversal walk} has transition matrix 
\[U=R(2\widehat{N}\widehat{N}^T-I),\]
where $R$ is the involution that reverses each arc, and $\widehat{N}$ is the normalized arc-tail incidence matrix. This model was formally introduced by Kendon \cite{Kendon2003}, and has been extensively studied in search algorithms. The famous Grover's search algorithm, for instance, is equivalent to an arc-reversal walk on the looped complete graph. 

We remark that each step of a vertex-face walk is equivalent to two steps of the arc-reversal walk, one on the original graph and one on the dual graph. Our computation shows that search using a vertex-face walk has a higher success probability than search using the arc-reversal walk. Of course, one reason is that the walker may move to non-adjacent arcs during each iteration of the vertex-face walk.

\section{Future Work}
The definition of vertex-face walks can be extended to non-orientable embeddings. To do this, we introduce graph-encoded maps. Let $\cM$ be a circular embedding, not necessarily orientable. A \textsl{flag} is an incident triple $(u,e,f)$ of vertex $u$, edge $e$ and face $f$. Pictorially, a flag is a triangle in the barycentric division of a face.  Figure \ref{Fig_flagsC3} gives the planar embedding of $C_3$, for which every dot represents a flag.

\begin{figure}[h]
	\centering
\begin{minipage}[b]{0.5\textwidth}
\begin{tikzpicture}
[every node/.style={circle,draw}]

\node (1) at (0,3.2){};
\node (2) at (-2,0){};
\node (3) at (2,0){};

\tikzset{inv/.style={inner sep=0pt}}

\node[inv] (0) at (0,1){};
\node[inv] (1') at (0,0){};
\node[inv] (2') at (1,1.6){};
\node[inv] (3') at (-1,1.6){};
\node[inv] (a) at (0,5.2){};
\node[inv] (a') at (0,-1.6){};
\node[inv] (b) at (-4.8,-1.6){};
\node[inv] (b') at (3,2.7){};
\node[inv] (c) at (4.8,-1.6){};
\node[inv] (c') at (-3,2.7){};

\tikzset{dot/.style={draw, inner sep=1.2pt, fill}}

\node[dot] (A) at (-0.3,2){};
\node[dot] (B) at (0.3,2){};
\node[dot] (C) at (1.1,0.8){};
\node[dot] (D) at (0.8,0.3){};
\node[dot] (E) at (-0.8,0.3){};
\node[dot] (F) at (-1.1,0.8){};

\node[dot] (G) at (-0.3,4){};
\node[dot] (H) at (0.3,4){};
\node[dot] (I) at (3,-0.2){};
\node[dot] (J) at (2.7,-0.6){};
\node[dot] (K) at (-2.7,-0.6){};
\node[dot] (L) at (-3,-0.2){};

\foreach \a/\b in {1/2, 2/3, 3/1}
\draw (\a) to (\b);

\foreach \a/\b in {0/1, 0/1', 0/2, 0/2', 0/3, 0/3', 1/a,1'/a',2'/b', 3'/c',2/b,3/c}
\draw [dashed] (\a) to (\b);
\end{tikzpicture}
\captionof{figure}{Planar embedding of $C_3$ and the flags}
\label{Fig_flagsC3}
\end{minipage}%
\begin{minipage}[b]{0.5\textwidth}
\begin{tikzpicture}
[every node/.style={circle,draw}]

\node (1) at (0,3.2){};
\node (2) at (-2,0){};
\node (3) at (2,0){};

\tikzset{inv/.style={inner sep=0pt}}

\node[inv] (0) at (0,1){};
\node[inv] (1') at (0,0){};
\node[inv] (2') at (1,1.6){};
\node[inv] (3') at (-1,1.6){};
\node[inv] (a) at (0,5.2){};
\node[inv] (a') at (0,-1.6){};
\node[inv] (b) at (-4.8,-1.6){};
\node[inv] (b') at (3,2.7){};
\node[inv] (c) at (4.8,-1.6){};
\node[inv] (c') at (-3,2.7){};

\tikzset{dot/.style={draw, inner sep=1.2pt, fill}}

\node[dot] (A) at (-0.3,2){};
\node[dot] (B) at (0.3,2){};
\node[dot] (C) at (1.1,0.8){};
\node[dot] (D) at (0.8,0.3){};
\node[dot] (E) at (-0.8,0.3){};
\node[dot] (F) at (-1.1,0.8){};

\node[dot] (G) at (-0.3,4){};
\node[dot] (H) at (0.3,4){};
\node[dot] (I) at (3,-0.2){};
\node[dot] (J) at (2.7,-0.6){};
\node[dot] (K) at (-2.7,-0.6){};
\node[dot] (L) at (-3,-0.2){};

\foreach \a/\b in {1/2, 2/3, 3/1}
\draw (\a) to (\b);

\foreach \a/\b in {0/1, 0/1', 0/2, 0/2', 0/3, 0/3', 1/a,1'/a',2'/b', 3'/c',2/b,3/c}
\draw [dashed] (\a) to (\b);

\foreach \a/\b in{B/C, D/E, F/A, H/I, J/K, L/G}
\draw[red, thick] (\a) to (\b);

\foreach \a/\b in {A/B, C/D, E/F, G/H, I/J, K/L}
\draw[blue, thick] (\a) to (\b);

\foreach \a/\b in {A/G, B/H, C/I, D/J, E/K, F/L}
\draw[green!70!black, thick] (\a) to (\b);
\end{tikzpicture}
\captionof{figure}{Planar embedding of $C_3$ and the gem}
\label{Fig_gemC3}
\end{minipage}
\end{figure}

For each flag $(u,e,f)$, let $u'$ be the other endpoint of $e$, let $e'$ be the other edge in $f$ that is incident to $u$, and let $f'$ be the other face that contains $e$. Define three functions
\begin{align*}
\tau_0: &(u, e, f) \mapsto (u', e, f),\\
\tau_1: &(u, e, f) \mapsto (u, e', f),\\
\tau_2: &(u, e, f) \mapsto (u, e, f').
\end{align*}
We have the following observations.
\begin{enumerate}[(i)]
\item $\tau_0, \tau_1, \tau_2$ are fixed-point-free involutions.
\item $\tau_0\tau_2 = \tau_2\tau_0$, and $\tau_0\tau_2 $ is fixed-point-free.
\item The group $\grp{\tau_0, \tau_1, \tau_2}$ acts transitively on the flags.
\end{enumerate}

If we join two flags in Figure \ref{Fig_flagsC3} by an edge whenever they are swapped by one of $\tau_0$, $\tau_1$ and $\tau_2$, then we obtain a cubic graph with a $3$-edge-coloring, as shown in Figure \ref{Fig_gemC3}.

In general, for a circular embedding $\cM$, a \textsl{graph-encoded map}, or \textsl{gem}, is a cubic graph with a $3$-edge coloring, where the vertices are the flags, and the $3$-edge coloring is induced by the three involutions $\tau_0$, $\tau_1$ and $\tau_2$, as described above. The concept of gem was first introduced by Lins in \cite{Lins1982}, where he also proved the following characterization of orientability.

\begin{theorem}
An embedding is orientable if and only if the gem is bipartite.
\qed
\end{theorem}

Note that an embedding $\cM$ with $\ell$ edges has $4\ell$ flags. Thus, if $\cM$ is orientable, then there are two components in the distance-$2$ graph of the gem, each with $2\ell$ vertices. Let $Y$ be one such component. We claim that the vertex-face walk for $\cM$ is equivalent to a quantum walk on the vertices of $Y$. Let $\pi_1$ be the partition of the vertices $(u,e,f)$ of $Y$ based on their third coordinates $f$. It is not hard to see that the size of each cell in $\pi_1$ is the degree of some face. Similarly, let $\pi_2$ be the partition of $V(Y)$ based on their first coordinates $u$. Let $\widehat{M}$ and $\widehat{N}$ be the normalized characteristic matrices for $\pi_1$ and $\pi_2$, respectively. Then 
\[(2\widehat{M}\widehat{M}^T-I)(2\widehat{N}\widehat{N}^T-I)\]
is precisely the vertex-face walk for $\cM$ relative to one consistent orientation of the faces. 

In general, let $\pi_1$ be coarsest partition of the flags for some circular embedding $\cM$, such that in each cell, all flags share an face, while no two flags share an edge. Similarly, let $\pi_2$ be the coarsest partition of the flags, such that in each cell, all flags share a vertex, while no two flags share an edge. Let $\widehat{M}$ and $\widehat{N}$ be the normalized characteristic matrices for $\pi_1$ and $\pi_2$, respectively. Then 
\[U=(2\widehat{M}\widehat{M}^T-I)(2\widehat{N}\widehat{N}^T-I)\]
defines a quantum walk, which is reducible if and only if $\cM$ is orientable. Now, each arc $(u,v)$ in the underlying graph $X$ is paired with two flags $(u,e,f)$ and $(u,e,f')$, and the probability that the walker is on the arc $(u,v)$ can be computed by summing the probabilities of her being on $(u,e,f)$ and $(u,e,f')$. There are many questions we may ask about this new definition of vertex-face walks; for example, one may study the relation between the original graph $X$ and the $H$-digraph, or compare the dynamics of vertex-face walks between orientable and non-orientable embeddings.

\bibliographystyle{amsplain}
\bibliography{dqw.bib}
\end{document}